\documentclass{article}
\usepackage{amsmath,amssymb,mathrsfs,accanthis,newpxmath}

\title{Limit distribution of errors in discretization of stochastic Volterra equations with multidimensional kernel}

\author{%
Masaaki Fukasawa
and
  Minato Hojo\\
{\small  Graduate School of Engineering Science, The University of Osaka}}

\date{}
\newcommand{\bm}[1]{\mathbf{#1}}
\newtheorem{theorem}{Theorem}
\newtheorem{lemma}{Lemma}
\newtheorem{remark}{Remark}

\begin{document}

\maketitle

\begin{abstract}
This paper investigates the limit distribution of discretization errors in stochastic Volterra equations (SVEs) with general multidimensional kernel structures. While prior studies, such as Fukasawa and Ugai (2023), were focused on one-dimensional fractional kernels, this research generalizes to broader classes, accommodating diagonal matrix kernels that include forms beyond fractional type. The main result demonstrates the stable convergence in law for the rescaled discretization error process, and the limit process is characterized under relaxed assumptions. 
\end{abstract}

\section{Introduction}

Stochastic Volterra equations (SVEs) 
 \begin{align}\label{SVE}
        X_t = X_0 + 
        \int_0^t \phi(t - s)b(X_s) \, ds + 
        \int_0^t \phi(t - s)\sigma(X_s) \, dW_s
    \end{align}
generalize stochastic differential equations (SDEs) by incorporating a Volterra kernel $\phi(t-s)$, allowing past states to influence the present. This property makes SVEs particularly suitable for modeling non-Markovian behaviors seen in fields like finance, neuroscience, and engineering. A prominent example is in rough volatility models~\cite{RoughVolatility}, where SVEs capture anti-persistent volatility behaviors of asset prices.

The study of discretization errors for SDEs is well-established, with significant results on their limit distributions (e.g.,\cite{jacod_protter_1998}). For SVEs, attention has primarily been given to fractional kernels, as demonstrated by \cite{FukasawaUgai2023, NualartSaikia}, which analyzed one-dimensional fractional kernels \(\phi(u) = u^{H-1/2}/\Gamma(H+1/2)\) with \(H \in (0, 1/2]\). 

This paper extends the framework to more general kernel structures, specifically diagonal matrix kernels \(\phi = \mathrm{diag}(\phi_1, \dots, \phi_d)\) that include forms beyond fractional type. These kernels introduce greater modeling flexibility while preserving the essential non-Markovian nature of SVEs. 
In particular, a local-stochastic rough volatility model~\cite{RoughVolatility}
    \begin{equation*}
        \begin{split}
           & d X^1_t = \sigma^1_1(X_t)\,dW^1_t + 
            \sigma^1_2(X_t)\,dW^2_t,
            \\
            & X^2_t = X^2_0 +\int_0^t (t-s)^{H-1/2}b^2(X_s)\,ds
             +\int_0^t (t-s)^{H-1/2}\sigma^2_1(X_s)\,dW^1_s
        \end{split}
    \end{equation*}
    falls into this generalized framework with
    $\phi(t) = \mathrm{diag}(1,t^{H-1/2})$.
The contribution of this study is to  establish the stable convergence in law for the rescaled discretization error process \( U^n = n^H(X - \hat{X}) \), where
    \begin{align}\label{approximation}
        \hat{X}_t = X_0 + 
        \int_0^t \phi(t - s)b\left(\hat{X}_{\frac{[ns]}{n}}\right) \, ds + 
        \int_0^t \phi(t - s)\sigma\left(\hat{X}_{\frac{[ns]}{n}}\right) \, dW_s
    \end{align}
    and $H \in (0,1)$ is determined by $\phi$, extending prior results to a generalized kernel framework with relaxed assumptions. 
By unifying and extending the approaches of earlier works, this paper lays a foundation for broader applications of SVEs in complex systems with non-Markovian dynamics.

\section{Main Result}\label{main}
Let $(\Omega, \mathscr{F}, \mathsf{P}, \{\mathscr{F}_t\}_{t \geq 0})$ be a filtered probability space satisfying the usual conditions.
Let $W$ is an $m$-dimensional standard Brownian motion defined on this space and assume that
    $X$ and  $\hat{X}$
   respectively   satisfy equations (\ref{SVE}) and (\ref{approximation}) for
     $b : \mathbb{R}^d \to \mathbb{R}^d$ and $\sigma : \mathbb{R}^d \to \mathbb{R}^{d \times m}$
            and $\phi : \mathbb{R} \to \mathbb{R}^{d\times d}$.
 We assume that the functions $b$ and $\sigma$ are continuously differentiable, with bounded and uniformly continuous derivatives.
We use the following notation:
\begin{itemize}
    \item \( C_0 \): The set of \( \mathbb{R}^d \)-valued continuous functions on \([0, T]\) vanishing at \( t = 0 \).
    \item \( C^\lambda_0 \): The set of \( \mathbb{R}^d \)-valued \( \lambda \)-H\"older continuous functions on \([0, T]\) vanishing at \( t = 0 \).
    \item \( \|\cdot\|_\infty \): The supremum norm on \([0, T]\).
    \item \( \|\cdot\|_{C^\lambda_0} \): The H\"older norm on \([0, T]\).
    \item \( \|\cdot\|_{L^p} \): The \( L^p \) norm with respect to \( P \).
    \item For any matrix \( A \), \( A^\top \) denotes the transpose of \( A \).
\end{itemize}
We introduce the following condition on the diagonal kernel $\phi=\mathrm{diag}(\phi_1,\phi_2,\cdots,\phi_d)$
for $H \in (0,1)$,
$\alpha\in[(1/2-H)\vee0,1/2)$ and $c_i \in \mathbb{R}$, $i=1,\dots,d$.\\

\noindent
{\bf A-$(H,\alpha,c_1,\dots,c_d)$:}
    There exist $\hat{H} \in (H,1)$ and $\hat{\phi}_i:(0,T] \to \mathbb{R}$, $i=1,\dots,d$ 
    such that
    \begin{itemize}
    \item  $\phi_i(u)=c_i u^{H-1/2}+\hat{\phi}_i(u)$, 
    \item $\hat{\phi}_i$ is absolutely continuous,
    \item$\hat{\phi}_i(u) = O(u^{\hat{H}-1/2})$ as $u \to 0$,
    \item$\hat{\phi}^{\prime}_i(u) = O(u^{\hat{H}-3/2})$ as $u \to 0$,
    \item$\mathcal{J}_{\phi_i}$ is continuous from $C^\lambda_0$ to $C^{\lambda-\alpha}_0$ for any $\lambda \in (\alpha,1/2)$
    \end{itemize}
    for each $i=1,\dots,d$,
    where 
$$\mathcal{J}_{\phi_i}f(t)=
\int_0^t \phi_i(t-s)df(s):=
\phi_i(t)f(t) - \int_0^t \phi_i^\prime(t - s)(f(t) - f(s)) \, ds.$$
Note that for any continuous It\^o process $Y$ with $Y_0 = 0$,
 \[
     (\mathcal{J}_{\phi_i} Y)(t) = \int_0^t \phi_i(t - s) \, dY_s 
    \]
    for all \( t \in [0, T] \) almost surely; see Proposition~A.2 of \cite{FukasawaUgai2023}.
A sufficient condition for {\bf A-$(H,\alpha, c_1,\dots,c_d)$} to hold with
$\alpha = (1/2-H) \vee 0$ and 
\begin{equation}\label{ci}
    c_i = \lim_{u \downarrow 0} u^{1/2-H}\phi_i(u)
\end{equation}
is that
$u \mapsto u^{1/2-H}\phi_i(u)$ is Lipschitz continuous for each $i=1,\dots, d$, as shown in Lemma~\ref{lemB} later.
The main result of this study is summarized in the following theorem:

\begin{theorem}\label{main result}
Assume {\bf A-$(H,\alpha,c_1,\dots,c_d)$} to hold and
    let  $\epsilon \in (0 , 1/2-\alpha)$. The process $\bm{U}^n = n^H (X - \hat{X})$ stably converges in law in $C_0^{1/2-\alpha-\epsilon}$ to a process $U = (U^1, \ldots, U^d)$,
    which is the unique continuous solution of the SVE
    \begin{equation}\label{mainu}\begin{split}
        U^i_t &= \sum_{k=1}^d \int_0^t \phi_i(t-s) U^k_s \left( \partial_k b^i(X_s) \, ds + \sum_{j=1}^m \partial_k \sigma^i_{j}(X_s) \, dW^j_s \right)  \\
        &-\frac{\Gamma(H+1/2)}{\sqrt{\Gamma(2H + 2) \sin (\pi H)}}\sum_{k=1}^d c_k \sum_{j=1}^m  \sum_{l=1}^m \int_0^t \phi_i(t-s) \partial_k \sigma^i_{j}(X_s) \sigma^k_{l}(X_s) \, dB^{l,j}_s,
    \end{split}\end{equation}
    where $B$ is an $m^2$-dimensional standard Brownian motion independent of $\mathscr{F}$ defined on 
    some extension of $(\Omega, \mathscr{F}, \mathsf{P})$.
\end{theorem}
{\it Proof. }
As in \cite{FukasawaUgai2023},
we set $\bm{U}^n=(U^{n,1},U^{n,2},\cdots,U^{n,d})$,
\[
V^{n,k,j} = n^H \int_0^{\cdot} (\hat{X}^k_s - \hat{X}^k_{\frac{[ns]}{n}}) \, dW^j_s
\]
for $1\leq k \leq d$, $1\leq j\leq m$
and $\Delta^n=(\Delta^{n,1},\cdots,\Delta^{n,d})^{T}$ by
\begin{align*}
    \Delta^{n,i}_{t}= &\, U^{n,i}_{t}  
    -\int_0^t \phi_i(t - s)
    \left(
    \nabla b^i\big(\hat{X}_{s}\big)^\top \bm{U}^{n}_{s} \, ds
    + \sum_{j=1}^m \nabla \sigma^{i}_{j}\big(\hat{X}_{s}\big)^\top \bm{U}^{n}_s \, dW^j_s
    \right) \\
    &-
    \int_0^t \phi_i(t - s) n^H \nabla b^i\big(\hat{X}_{s}\big)^\top 
    \big(\hat{X}_s - \hat{X}_{\frac{[ns]}{n}}\big) \, ds \\
   & -
    \sum_{j=1}^m \sum_{k=1}^d \int_0^t \phi_i(t - s) \partial_k \sigma^i_{j}\big(\hat{X}_{s}\big) \, dV^{n,k,j}_s.
\end{align*}
    The result then follows by combining Lemmas~2.2-7 below
    as detailed in~\cite{FukasawaUgai2023}.
\hfill{$\square$}
  \begin{lemma}\label{2.3}
    \( \bm{V}^n :=\{ V^{n,k,j} \}_{1\leq k\leq d,1\leq j\leq m}\) stably converges in law in \( C_0 \) and 
    the limit \( V = \{ V^{k,j} \} \) can be expressed as
    \[
    V^{k,j} =\frac{\Gamma(H+1/2)}{\sqrt{\Gamma(2H + 2) \sin (\pi H)}}c_k \sum_{l=1}^{m}  \int_0^t \sigma^k_l(X_s) \, dB_s^{l,j}
    \]
    where \( B \) is an \( m^2 \)-dimensional standard Brownian motion, independent of 
    \( \mathscr{F} \) and defined on some extension of \( (\Omega, \mathscr{F}, \mathsf{P}) \).
    \end{lemma}

    \begin{lemma}\label{2.4}
        For all \( i \in \{1, \ldots, d\} \), for any $\epsilon\in(0,1/2-\alpha)$
        \[
        \int_0^t \phi(t-s) n^H \nabla b^i(X_{s})^\top (X_s - \tilde{X}_s) \, ds 
        \to 0 \quad \text{in probability, in } C_0^{1/2-\alpha-\epsilon}.
        \]
    \end{lemma}
        
    \begin{lemma}\label{2.5}
        $\|\Delta^n\|_{C_0^{\gamma}} $ tends to zero in \( L^p \) for any \( \gamma \in (0, 1/2-\alpha) \) and \( p \geq 1 \).
    \end{lemma}
         
    \begin{lemma}\label{2.6}
        If the sequence
        \[
        (\bm{U}^n, \bm{V}^n, \{\nabla b^i(\hat{X})\}_i, \{\nabla \sigma^i_j(\hat{X})\}_{ij})
        \]
        converges in law in 
        $
        C_0^{1/2-\alpha-\epsilon} \times C_0 \times (C_0)^d \times (C_0)^{dm}
        $
        to
        \[
        (U, V, \{\nabla b^i(X)\}_i, \{\nabla\sigma^i_j(X)\}_{ij}),
        \]
        then \( U \) is the solution of (2.1).
    \end{lemma}
        
    \begin{lemma}\label{2.7}
        The sequence \( \bm{U}^n \) is tight in \( C_0^{H-\epsilon} \) for any \( \epsilon \in (0, H) \).
    \end{lemma}
        
    \begin{lemma}\label{2.8}
        The uniqueness in law holds for continuous solution of (\ref{mainu}).
    \end{lemma}

 The proofs of these lemmas are omitted because they are straightforward extensions of Lemmas~2.3-8 of \cite{FukasawaUgai2023} after Lemmas~\ref{2.2} and 3.1-7 below are established.

  \begin{lemma}\label{2.2}
    For all $t \in [0, T]$, $(k_1, k_2) \in \{1, \ldots, d\}^2$, and $1 \leq j \leq m$,
        \begin{enumerate}
            \item[(i)] 
            \[
            \langle V^{n,k_1,j}, V^{n,k_2,j} \rangle_t 
            \xrightarrow{L^1} 
           \frac{\Gamma(H+1/2)^2}{\Gamma(2H + 2) \sin (\pi H)} c_{k_1}c_{k_2}
            \sum_{l=1}^m \int_0^t \sigma^{k_1}_l(X_s) \sigma^{k_2}_l(X_s) \, ds,
            \]
            \item[(ii)] 
            \[
            \langle V^{n,k,j}, W^j \rangle_t 
            \xrightarrow{L^1} 0, 
            \]
        \end{enumerate}
        as $n\to \infty$.
    \end{lemma}
The proof of this lemma is 
is given in Section 4.

\begin{remark}
The H\"older spaces are not separable.
However, according to Section 2.1 of~\cite{hamadouche_2000},  the $\beta$-H\"older space can be regarded as a separable subspace of the $\gamma$-H\"older space for $\gamma<\beta$.
This property resolves all the delicate measurablity issues for nonseparable-space-valued random variables in this study.
\end{remark}

\section{Preliminary estimates}
\subsection{Estimates for the kernel}
Here we derive a few estimates which play a key role in this study. We set $\bar{H}=H/2+1/2$, $\beta \in (1,(1-2H)^{-1})$ for $H \in (0,1/2)$, $\beta=2$ for $H \in [1/2,1)$ and $\beta^* = \beta/(\beta - 1)$.
These satisfy $\bar{H}>H$, $\int_0^t |\phi_i(s)|^{2\beta} \, ds<\infty$ and $1/\beta+1/\beta^*=1$.
We use $C$ to represent a constant which may differ line by line.
\begin{lemma}\label{orders of phi}
There exists $\bar{H} \in (H,1)$ such that
    \begin{align*}
        \int_0^h |\phi_i(t)| \, dt = O(h^{H+1/2}),
        \quad
        \int_0^T |\phi_i(t + h) - \phi_i(t)| \, dt = O(h^{\bar{H}}),\\
        \left( \int_0^h |\phi_i(t)|^2 \, dt \right)^{1/2} = O(h^H),
        \quad
        \left( \int_0^T |\phi_i(t + h) - \phi_i(t)|^2 \, dt \right)^{1/2} = O(h^H).
    \end{align*}
\end{lemma}
{\it Proof. }
    By the {\bf A-$(H,\alpha,c_1,\dots,c_d)$}, $|\phi(t)|\leq C t^{H-1/2}$.
    This leads to $\int_0^h |\phi_i(t)| \, dt = O(h^{H+1/2})$ and $\left( \int_0^h |\phi_i(t)|^2 \, dt \right)^{1/2} = O(h^H)$.
    In addition, we can see that for any $\alpha \in (0,H]$,
    \begin{align*}
        |\phi_i(t + h) - \phi_i(t)|&=|\int_t^{t+h} \phi^{\prime}(u)\,du|
        \leq \int_t^{t+h} |\phi^{\prime}(u)|\,du\\
        &\leq \int_t^{t+h} u^{\alpha-3/2} \,du \leq C|(t+h)^{\alpha-1/2}-t^{\alpha-1/2}|,
    \end{align*}
    so we conclude that by change of variables $u=ht$
    \begin{align*}
        \int_0^T |\phi_i(t + h) - \phi_i(t)| \, dt
        &\leq C\int_0^T |(t + h)^{H/2-1/2} - t^{H/2-1/2}| \, dt\\
        &\leq C h^{\bar{H}} \int_0^{\infty}|(u+1)^{H/2-1/2} -u^{H/2-1/2}|\,du  = O(h^{\bar{H}}).
    \end{align*}
    Similarly, we have
    \begin{align*}
        \int_0^T (\phi_i(t + h) - \phi_i(t))^2 \, dt
        &\leq C h^{2H} \int_0^{\infty}((u+1)^{H-1/2} -u^{H-1/2})^2\,du = O(h^{2H})
    \end{align*}
    which concludes the proof.
\hfill{$\square$}
The following lemma is proven in the same way, so the proof is omitted.
\begin{lemma}\label{orders of phi hat}
    \begin{align*}
        \left( \int_0^h |\hat{\phi}_i(t)|^2 \, dt \right)^{1/2} = O(h^{\hat{H}}),
        \quad
        \left( \int_0^T |\hat{\phi}_i(t + h) - \hat{\phi}_i(t)|^2 \, dt \right)^{1/2} = O(h^{\hat{H}}).
    \end{align*}
\end{lemma}
The following lemma is derived from
Lemmas~\ref{orders of phi} and \ref{orders of phi hat}
in the same manner as Lemma 3.1 of \cite{FukasawaUgai2023}, 
so the proof is omitted.
\begin{lemma}\label{Kernela}
    The following inequalities hold for any adapted \( \mathbb{R}^d \)-valued process \( Y \) and \( \mathbb{R}^{d \times m} \)-valued process \( Z \):
    \begin{enumerate}
        \item For \( p \geq 2 \),
        \[
        \mathbb{E} \left[ \left| \int_0^t \phi(t-s) Y_s \, ds \right|^p \right]
        \leq C \int_0^t \mathbb{E} \left[ |Y_s|^p \right] \, ds,
        \]
    
        \item For \( p > 2 \beta^* \),
        \[
        \mathbb{E} \left[ \left| \int_0^t \phi(t-s) Z_s \, dW_s \right|^p \right]
        \leq C \int_0^t \mathbb{E} \left[ |Z_s|^p \right] \, ds,
        \]
        \item For \( p \geq 1 \),
        \begin{align*}
        \mathbb{E} &\left[ \left| \int_0^t (\phi(t+h-s) - \phi(t-s)) Y_s \, ds \right|^p \right]
        + \mathbb{E} \left[ \left| \int_t^{t+h} \phi(t+h-s) Y_s \, ds \right|^p \right]\\
        &\leq C h^{\bar{H}p} \sup_{r \in [0,T]} \mathbb{E} \left[ |Y_r|^p \right].
        \end{align*}
    
        \item For \( p \geq 2 \),
        \begin{align*}
        \mathbb{E}& \left[ \left| \int_0^t (\phi(t+h-s) - \phi(t-s)) Z_s \, dW_s \right|^p \right]
        + \mathbb{E} \left[ \left| \int_t^{t+h} \phi(t+h-s) Z_s \, dW_s \right|^p \right]\\
        &\leq C h^{Hp} \sup_{r \in [0,T]} \mathbb{E} \left[ |Z_r|^p \right].
        \end{align*}
    \end{enumerate}
    Here, the constant \( C \) depends only on \( K, \beta, p, \) and \( T \).
\end{lemma}

\subsection{Intermediate results}
The following lemmas are presented as Lemmas 3.5-8 in \cite{FukasawaUgai2023} under different conditions on the kernel. By using Lemma \ref{Kernela}, they can be proven in the same way, so their proofs are omitted.
\begin{lemma}\label{sup x hat bdd}
    Let $p \geq 1$. Then,
    \[
    \sup_{t \in [0, T]} E\left[ |\hat{X}_t|^p \right] \leq C,
    \]
    where $C$ is a constant that only depends on $|X_0|$, $|b(0)|$, $|\sigma(0)|$, $K$, $p$, $\beta$, and $T$.
\end{lemma}

\begin{lemma}\label{Lipschitz of X hat}
    Let $p \geq 1$. Then,
    \[
    E\left[ |\hat{X}_t - \hat{X}_s|^p \right] \leq C |t - s|^{Hp}, \quad t, s \in [0, T],
    \]
    and for $p > H^{-1}$,
    \[
    E \left[ \sup_{0 \leq s \leq t \leq T} \frac{|\hat{X}_t - \hat{X}_s|}{|t - s|^\gamma}^p \right] \leq C_\gamma
    \]
    for all $\gamma \in [0, H - p^{-1})$, where $C_\gamma$ is a constant that does not depend on $n$. As a consequence, $\hat{X}$ is a $C^\gamma$-valued random variable for any order $\gamma < H$ for all $n$.
\end{lemma}

\begin{lemma}\label{difference of x and x hat}
    Let $p \geq 1$. Then the process $X_t - \hat{X}_t$ uniformly converges to zero in $L^p$ with the rate $n^{-Hp}$ as $n$ goes to infinity, that is,
    \[
    \sup_{t \in [0, T]} E\left[ |X_t - \hat{X}_t|^p \right] \leq C n^{-Hp},
    \]
    where $C$ is a positive constant which does not depend on $n$.
\end{lemma}

\begin{lemma}
    For all $p \geq 1$ and $\epsilon \in (0, H)$, there exists a constant $C > 0$ which does not depend on $n$ such that
    \[
    E \left[ \sup_{t \in [0, T]} |X_t - \hat{X}_t|^p \right] \leq C n^{-p(H - \epsilon)}.
    \]
\end{lemma}
    
    \section{Proof of Lemma~\ref{2.2}}

    We first introduce the following definitions
\begin{align*}
    &\psi^{n,k}_{1,s}:=\int_{0}^{\frac{[ns]}{n}} (\phi_k(s-u)-\phi_k(\frac{[ns]}{n}-u))b^k(\tilde{X}_u) \, du,\\
    &\psi^{n,k}_{2,s}:=b^k(\tilde{X}_s) \int_{\frac{[ns]}{n}}^{s} \phi_k(s-u) \, du ,\\
    &\psi^{n,k,j}_{3,s}:=c_k\int_{0}^{\frac{[ns]}{n}} ((s-u)^{H-1/2}-(\frac{[ns]}{n}-u)^{H-1/2}) {\sigma}^k_{j}(\tilde{X}_u) \, dW^j_u ,\\
    &\psi^{n,k,j}_{4,s}:=c_k{\sigma}^k_{j}(\tilde{X}_s) \int_{\frac{[ns]}{n}}^{s} (s-u)^{H-1/2} \, dW^j_u,\\
    &\psi^{n,k,j}_{5,s}:=\int_{0}^{\frac{[ns]}{n}} (\hat{\phi}_k(s-u)-\hat{\phi}_k(\frac{[ns]}{n}-u)) \sigma^k_{j}(\tilde{X}_u) \, dW^j_u ,\\
    &\psi^{n,k,j}_{6,s}:=\sigma^k_{j}(\tilde{X}_s) \int_{\frac{[ns]}{n}}^{s} \hat{\phi}_k(s-u) \, dW^j_u
\end{align*}
for
$k=1,2,\cdots,d $ and $j=1,2,\cdots,m$,
where $\tilde{X}_t = \hat{X}_{[nt]/n}$.
Observe that
\[
    \hat{X}^k_s - \hat{X}^k_{\frac{[ns]}{n}}=\psi^{n,k}_{1,s}+\psi^{n,k}_{2,s}
    +\sum_{j=1}^m \psi^{n,k,j}_{3,s}+\sum_{j=1}^m \psi^{n,k,j}_{4,s}
    +\sum_{j=1}^m \psi^{n,k,j}_{5,s}+\sum_{j=1}^m \psi^{n,k,j}_{6,s}.
\]
By Lemma~4.2 of \cite{FukasawaUgai2023}, we have
\begin{align*}
      &\sup_{n\geq0}\sup_{s\in[0,T]}n^H\|\psi^{n,k,j}_{3,s} +\psi^{n,k,j}_{4,s} \|_{L^2}<\infty 
\end{align*}
and
  \begin{align*}
          n^{2H}&\sum_{j,l=1}^m\int_0^t(\psi^{n,k_1,j}_{3,s}+ \psi^{n,k_1,j}_{4,s})(\psi^{n,k_2,l}_{3,s}+ \psi^{n,k_2,l}_{4,s})\, ds\\ &\to 
       \frac{\Gamma(H+1/2)^2}{\Gamma(2H + 2) \sin (\pi H)} c_{k_1}c_{k_2} \delta^{jl}
             \int_0^t {\sigma}_{k_1}^j(X_s) {\sigma}_{k_2}^j(X_s) \, ds
    \end{align*}
in $L^1$ as $n\to \infty$, where $\delta^{jl}$ is the Kronecker 
delta.
Note that $H\leq 1/2$ is assumed in \cite{FukasawaUgai2023} and used only in
Lemma~4.1 of \cite{FukasawaUgai2023}.
To include the case $H>1/2$,
we provide Lemmas~\ref{A} and \ref{B} below.

Now, in order to prove Lemma~\ref{2.2}-(i), it suffices then to show the following lemma.
\begin{lemma}  \label{lim of phi}
    For $k=1,2,\cdots,d$ and $j=1,2,\cdots,m$
    \begin{align*}
        &\lim_{n\to\infty}\sup_{s\in[0,T]}  n^H\|\psi^{n,k}_{i,s}\|_{L^2}=0 \quad i=1,2,\\
        &\lim_{n\to\infty}\sup_{s\in[0,T]}n^H\|\psi^{n,k,j}_{i,s}\|_{L^2}=0 \quad i=5,6.  
    \end{align*}
\end{lemma}

{\it Proof. }
    By Minkowski's integral inequality, Lemmas~\ref{orders of phi} and \ref{sup x hat bdd} and change of variable, 
    \begin{align*}    
        \|\psi^{n,k}_{1,s}\|^2_{L^2} 
        &\leq \left( \int_0^{\frac{[ns]}{n}} \left|\phi_k(s-u) - \phi_k\left(\frac{[ns]}{n} - u\right)\right| E [ |b_k(\hat{X}_{\frac{[nu]}{n}})|^2 ]^{\frac{1}{2}} \, du \right)^2 \\
        &\leq \sup_{r \in [0, T]} E [ |b_k(\hat{X}_r)|^2 ] \left( \int_0^{\frac{[ns]}{n}} \left|\phi_k(u+s-\frac{[ns]}{n}) - \phi_k(u)\right| \, du \right)^2  \leq Cn^{-2\bar{H}}
        \end{align*}
        and
      \begin{align*}    
        \|\psi^{n,k}_{2,s}\|^2_{L^2}
        &\leq E \left[ \left( b_k (\hat{X}_{\frac{[ns]}{n}}) \int_{\frac{[ns]}{n}}^{s} \phi_k(s-u) \, du \right)^2 \right] \\
        &\leq \sup_{r \in [0, T]} E [ |b_k(\hat{X}_r)|^2 ] \left(\int_{\frac{[ns]}{n}}^{s}|\phi_k(s-u)| \, du\right)^2
        \\ & \leq C\left(\int^{s-\frac{[ns]}{n}}_{0}|\phi_k(u)| \, du\right)^2 \leq Cn^{-2(H+1/2)}.
    \end{align*}
    The Burkholder-Davis-Gundy inequality leads to following  inequalities in a similar way:
    \begin{align*}
        \|\psi^{n,k,j}_{5,s}\|^2_{L_2}
        &\leq \sup_{r \in [0, T]} E [(\sigma^k_j(\hat{X}_{r}))^2 ]\int_0^{\frac{[ns]}{n}} \left(\hat{\phi}_k(s-u) - \hat{\phi}_k\left(\frac{[ns]}{n} - u\right)\right)^2  \,du \leq Cn^{-2\hat{H}},\\
        \|\psi^{n,k,j}_{6,s}\|^2_{L^2}
        &\leq\sup_{r \in [0, T]} E [(\sigma^k_j(\hat{X}_{r}))^2 ]\int_{\frac{[ns]}{n}}^{s} \hat{\phi}_k(s-u)^2\, du \leq Cn^{-2\hat{H}}.
    \end{align*}
    These inequalities imply the assertion.
\hfill{$\square$}

To prove Lemma~\ref{2.2}-(ii),
we set $\Delta \hat{X}_s = \hat{X}_s - \hat{X}_{[ns]/n}$. We have
\[
    \langle V^{n,k,i}, W^i \rangle_t = \int_0^t n^H \Delta \hat{X}_s \, ds
\]
and by Fubini's theorem,
\[
    E\left[ \left| \langle V^{n,k,i}, W_i \rangle_t \right|^2 \right] = 2 \int_0^t \int_0^s n^{2H} E\left[\Delta \hat{X}^k_s \Delta \hat{X}^k_v \right]\, dv \, ds .
\]
We will check inequalities and convergences to use the dominated convergence theorem.
By Lemma~\ref{Lipschitz of X hat} and the Cauchy-Schwarz inequality,
\begin{align*}
    |E[\Delta \hat{X}^k_s \Delta \hat{X}^k_v]| &\leq E\left[|\Delta \hat{X}^k_s|^2\right]^{\frac{1}{2}} E\left[|\Delta \hat{X}^k_v|^2\right]^{\frac{1}{2}} \\ &\leq C \left(\frac{s - [ns]}{n}\right)^H \left(\frac{v - [nv]}{n}\right)^H \leq C n^{-2H}.
\end{align*}
We next show 
$n^{2H} E[\Delta \hat{X}^k_s \Delta \hat{X}^k_v ]\to 0$. From Lemma~\ref{lim of phi}, we deduce that
\[
    n^{2H}(\Delta \hat{X}^k_s \Delta \hat{X}^k_v-\sum_{j=1,l=1}^m(\psi^{n,k,j}_{3,s}+\psi^{n,k,j}_{4,s})(\psi^{n,k,l}_{3,v}+\psi^{n,k,l}_{4,v})) \to 0\quad \text{in } L_1.
\]
The result then follows as in Section~4.2 of \cite{FukasawaUgai2023}
using Lemmas~\ref{A} and \ref{B} below.
\hfill $\square$

\appendix 
\renewcommand{\thesection}{\Alph{section}} 

\section{Auxiliary lemmas}

\begin{lemma}\label{preA}
    Let $\alpha<1$, $0<x<y$, $y'\leq x'$ and $0\leq x'<x$, then
    $|y^\alpha-x^\alpha|\leq|(y-y')^\alpha-(x-x')^\alpha|$
\end{lemma}
{\it Proof. }
    Let $f(s,t)=|(s+t)^\alpha-s^\alpha|$ $(s>0,t\geq0)$, then
    \[
        \frac{\partial}{\partial s}f(s,t)=|\alpha|((t + s)^{-1 + \alpha}- s^{-1 + \alpha} )\leq0
    \]
    \[
        \frac{\partial}{\partial t}f(s,t)=|\alpha| (t + s)^{-1 + \alpha}\geq0.
    \]
    Therefore, we have
    \[
        f(x,y-x)\leq f(x-x',y-x) \leq f(x-x',y-x+(x'-y')).
    \]
    Since $|y^\alpha-x^\alpha|=f(x,y-x)$ and $|(y-y')^\alpha-(x-x')^\alpha|=f(x-x',y-x+(x'-y'))$,
    this proof is completed.

\hfill{$\square$}
    
\begin{lemma}\label{A}
    Let $\alpha \in (-1/2,1/2)$ and
    \[
        A_n(v,s)=n^{2\alpha+1} \int^{\frac{[nv]}{n}}_{0} \left( (s - u)^{\alpha} - \left( \frac{[ns]}{n} - u \right)^{\alpha} \right) \left( (v - u)^{\alpha} - \left( \frac{[nv]}{n} - u \right)^{\alpha} \right) du
    \]for $v \leq s$.
    Then $\sup_{0 \leq v\leq s\leq M}|A_n(v,s)|<\infty$ and $\lim_{n\to \infty}A_n(v,s)=0$
\end{lemma}
        
{\it Proof. }
    It is clear that $A_n(v,s)\geq 0$ from the assumption. By change of variable $z=[nv]-nu$, we have
    \begin{align*}
&A_n(v,s)\\
&=\int_0^{[nv]} \left( \left( z + ns - [nv] \right)^{\alpha} - \left( z + [ns] - [nv] \right)^{\alpha} \right) \left( \left( z + nv - [nv] \right)^{\alpha} - z^{\alpha} \right) \, dz.
    \end{align*}
    In addition, by considering $(x',y')$ in Lemma \ref{preA} to be $(ns - [nv]-1,[ns] - [nv])$ and $(nv - [nv]-1,0)$, we obtain
    \[
        |\left( z + ns - [nv] \right)^{\alpha} - \left( z + [ns] - [nv] \right)^{\alpha}|\leq|(z+1)^\alpha-z^\alpha|,
    \]
    \[
        |\left( z + nv - [nv] \right)^{\alpha} - z^{\alpha}|\leq|(z+1)^\alpha-z^\alpha|.
    \]
    By combining these two inequalities, we have
    \begin{align*}
 &1_{[0,[nv]]}(z)\left| \left( z + ns - [nv] \right)^{\alpha} - \left( z + [ns] - [nv] \right)^{\alpha} \right| \left| \left( z + nv - [nv] \right)^{\alpha} - z^{\alpha} \right|
    \\ & \leq((z+1)^\alpha-z^\alpha)^2.
    \end{align*}
    Also, it follows that
    \begin{align}
        | (z + ns - [nv] )^{\alpha} - ( z + [ns] - [nv] )^{\alpha}|
        &= |\alpha|\int_{z + ns - [nv]}^{z + [ns] - [nv]}w^{\alpha-1}\,dw\nonumber\\
        &\leq |\alpha|\int_{z + ns - [nv]}^{z + [ns] - [nv]}(z + ns - [nv])^{\alpha-1}\,dw \label{A112}\\
        &\leq |\alpha|(z + ns - [nv])^{\alpha-1} \to 0\text{ as } n\to\infty.\nonumber
    \end{align}
    As a result, we get
    \[
        |A_n(v,s)|\leq\int_{0}^{\infty}((x + 1)^\alpha - x^\alpha)^2\,dx
    \]
    and by the dominated convergence theorem, the proof is completed.
\hfill{$\square$}

\begin{lemma}\label{B}
    Let $\alpha \in (-1/2,1/2)$ and
    \[
    B_n(v,s):=n^{2\alpha+1}\int_{[nv]/n}^{v} \left|(s-u)^\alpha - \left(\frac{[ns]}{n} - u\right)^\alpha\right| (v-u)^\alpha \, du
    \]
    for $nv\leq [ns] $. Then 
    $ \lim_{n\to\infty} B(v,s)=0 $
\end{lemma}

{\it Proof. }
    By change of variable $z=n(v-u)$, it is holds that 
    \[
    B_n(v,s)=
    \int_{0}^{nv-[nv]} | (z+ ns - nv )^\alpha - (z+ [ n s ] - nv)^\alpha | z^\alpha \, dz
    \]
    In addition, by Lemma \ref{preA}, we have
    \[
        1_{(0,nv-[nv])}| (z+ ns - nv  )^\alpha - (z+ [ n s ] - nv  )^\alpha | z^\alpha
        <1_{(0,1)}| (z+1)^\alpha - z^\alpha | z^\alpha
    \]
    and in the same way as (\ref{A112}), 
    $| (z+ ns - nv  )^\alpha - (z+ [ n s ] - nv  )^\alpha | \to 0$
    so the dominated convergence theorem leads to the assertion.
\hfill{$\square$}

    \begin{lemma}
    \label{lemB}
        Let $H\in(0,1)$ and $\alpha=(1/2-H)\vee 0$. If $u^{1/2-H}\phi_i(u)$ is Lipschitz continuous on $(0,T]$ for each $i=1,\dots, d$, then 
         the condition {\bf A-$(H,\alpha,c_1,\dots,c_d)$} holds with
         \eqref{ci}.
    \end{lemma}
{\it Proof. }
We follow the proof of Lemma F.3 of \cite{Ugai2022} with a slight extension.
    By the assumption, each $\phi_i$ is expressed by
    \[
        \phi_i(u)=f(0)u^{H-1/2}-(f(u)-f(0))u^{H-1/2}
    \]
    for a Lipschitz continuous function $f$.
    Also, the Lipschitz continuity leads to
    \[
        |(f(u)-f(0))u^{H-1/2}|\leq C u\cdot u^{H-1/2}=Cu^{H+1/2}
    \]
    and
    \begin{align*}
        |((f(u)-f(0))u^{H-1/2})'|
        &\leq|f'(u)u^{H-1/2}|+|(f(u)-f(0))u^{H-3/2}|\\
        &\leq \sup_{s\in[0,T]}|f'(s)|u^{H-1/2}+C u^{H-1/2}.
    \end{align*}
    Therefore, it is sufficient to check the continuity of $\mathcal{J}_{\phi_i}$.

    Let $f(u) = u^\alpha \phi_i(u)$, $d_1(u)= f(u)u^{-(\alpha+1)}$, $d_2(u)=f'(u)u^{-\alpha}$ and
    for $g \in C^\lambda_0$,
    \begin{align*}
    &\mathcal{M}(g):=fg, \quad \mathcal{D}g(t):=\frac{g(t)}{t^{\alpha}},\\
    &\mathcal{I}_ig(t):=\int_0^t d_i(t - s)(g(t) - g(s)) \, ds,\quad i=1,2.
    \end{align*} 
    Then we have
    \[
        \mathcal{J}_{\phi_i}=\mathcal{D}\mathcal{M}+\mathcal{I}_1+\mathcal{I}_2.
    \]
    We will prove the continuity of each operator. Let $\lambda\in(\alpha,1)$.

\subsubsection*{Proof of the continuity of $\mathcal{M}$ from $C^\lambda_0$ to $C^{\lambda}_0$}
Let $t,s \in [0,T]$. 
We have $\mathcal{M}g(0)=f(0)g(0)=0$ and
\begin{align*}
    |\mathcal{M}g(t)-\mathcal{M}g(s)|
    \leq |f(t)||g(t)-g(s)|+|g(s)||f(t)-f(s)|
    \leq C\|g\|_{C^\lambda_0}|t-s|^{\lambda},
\end{align*}
since $\sup_{t\in{[0,T]}}|f(t)|<\infty$ and $f$ is 
$\lambda$--H\"older continuous.
Therefore we obtain the continuity.

\subsubsection*{Proof of the continuity of $\mathcal{D}$ from $C^\lambda_0$ to $C^{\lambda-\alpha}_0$}
Let $t,s \in (0,T]$ and $t>s$. We have
    \[
        |\mathcal{D} g(t)|\leq \frac{g(t)-g(0)}{t^{\alpha}}
        \leq \|g\|_{C^{\lambda}_0} t^{\lambda-\alpha},
    \]
    so $\mathcal{D} f$ can be defined on $[0,T]$ and $\mathcal{D}g(0)=0$.
    Next we evaluate the difference.
    \[
        |\mathcal{D} g(t)-\mathcal{D} g(s)| 
        = \left|\frac{g(t)-g(s)}{t^\alpha}\right| 
        + |g(s)|\left|\frac{1}{t^\alpha}-\frac{1}{s^\alpha}\right|.
    \]
    We will show that the first term is bounded. $|g(t)-g(s)|\leq \|g\|_{C^{\lambda}_0}|t-s|^{\lambda}$ and $t^{\alpha}\geq |t-s|^{\alpha}$ lead to
    \[
        \left|\frac{g(t)-g(s)}{t^\alpha}\right| 
        \leq \|g\|_{C^{\lambda}_0}|t-s|^{\lambda-\alpha}
    .\]
    Next we will show that the second term is bounded. 
    We have 
    \[
        |g(s)|\leq \|g\|_{C^{\lambda}_0} s^{-\lambda},
    \]
    and there exists a constant $C>0$ such that
    \[
        \left|\frac{1}{t^\alpha}-\frac{1}{s^\alpha}\right|
        \leq C s^{\lambda}|t-s|^{\lambda-\alpha}
    \]
    by the following argument.

    In the case where $2s>t$, we have 
    \begin{align*}
        \left|\frac{1}{t^\alpha}-\frac{1}{s^\alpha}\right|
        &\leq \alpha \int_{s}^{t}x^{-\alpha-1}\,dx \leq \alpha \int_{s}^{t}s^{-\alpha-1}\,dx =\alpha s^{-\alpha-1}|t-s|\\
        &\leq \alpha s^{-\alpha-1} |2s-s|^{1-(\lambda-\alpha)}|t-s|^{\lambda-\alpha}=\alpha s^{-\lambda}|t-s|^{\lambda-\alpha},
    \end{align*}
    and in the other case $(2s<t)$, we have 
    \begin{align*}
        |t^{-\alpha}-s^{-\alpha}|
        \leq s^{-\alpha}
        \leq s^{-\lambda}|t-s|^{\lambda-\alpha}.
    \end{align*}
    Therefore we conclude that the second term is bounded by a constant multiple of $\|g\|_{C^{\lambda}_0}|t-s|^{\lambda-\alpha}$.
    This implies the continuity.

\subsubsection*{Proof of the continuity of $\mathcal{I}_1$ from $C^\lambda_0$ to $C^{\lambda-\alpha}_0$}
    Let $h\in(0,1)$ and $t\in(0,T)$ such that $t+h\leq T$.
    We have
    \[
    |\mathcal{I}_1g(t)|\leq \|g\|_{C_{0}^{\lambda}}\sup_{s\in(0,T]}|f(s)|\,\int_0^{t}s^{\lambda-\alpha-1}\,ds
    \]
    so $\mathcal{I}_1 f$ can be defined on $[0,T]$ and $\mathcal{I}_1 g(0)=0$.
    Next we evaluate the difference.
    By the change of variable, we have
    \[
        \mathcal{I}_1g(t)=\int_{0}^{t}(g(t)-g(t-s))d_1(s)\,ds,
    \]
    \[
        \mathcal{I}_1g(t+h)=\int_{-h}^{t}(g(t+h)-g(t-s))d_1(s+h)\,ds.
    \]
    These lead to the inequality:
    \begin{align*}
        |&\mathcal{I}_1g(t+h)-\mathcal{I}_1g(t)|\\
        &\leq \int_{0}^{t} |g(t)-g(t-s)||d_1(s+h)-d_1(h)|\,ds\\
        &+ \int_{0}^{t} |g(t+h)-g(t)||d_1(s+h)|\,ds
        + \int_{-h}^{0}|g(t+h)-g(t-s)||d_1(s+h)|\,ds\\
        &\leq \|g\|_{C_{0}^{\lambda}}(
        \int_{0}^{t} |f(s+h)||(s+h)^{-\alpha-1}-s^{-\alpha-1}|\,ds+\int_{0}^{t} s^{\lambda}|f(s+h)-f(s)|\,ds
        \\
        &\quad\quad\quad\quad+ \int_{0}^{t} h^{\lambda}|d_1(s+h)|\,ds
        + \int_{-h}^{0}(s+h)^{\lambda}|d_1(s+h)|\,ds
        )
        \\
        &\leq C\|g\|_{C_{0}^{\lambda}}
        (
        \int_{0}^{t} s^{\lambda}|(s+h)^{-\alpha-1}-s^{-\alpha-1}|\,ds
        +h\int_{0}^{t} s^{\lambda-\alpha-1}\,ds
        \\
        &\quad\quad\quad\quad+ \int_{0}^{t} h^{\lambda}(s+h)^{-\alpha-1}\,ds
        + \int_{-h}^{0}(s+h)^{\lambda-\alpha-1}\,ds
        ).
    \end{align*}
    This is bounded by a constant multiple of $\|g\|_{C_0^{\lambda}}h^{\lambda-\alpha}$ because
        \[
    \int^{0}_{-h} (s + h)^{\lambda - \alpha - 1} \, ds =C h^{\lambda - \alpha}, \quad
        \int_{0}^\infty (s + h)^{-\alpha - 1} \, ds= C h^{-\alpha},
    \]
    \begin{align*}
        \int_{0}^{t} s^{\lambda}|(s+h)^{-\alpha-1}-s^{-\alpha-1}|\,ds
        &\leq h^{\lambda-\alpha}\int_{0}^{t/h}s^{\lambda}|(s+1)^{-\alpha-1}-s^{-\alpha-1}|\,ds\\
        &\leq h^{\lambda-\alpha}\int_{0}^{\infty}s^{\lambda}|(s+1)^{-\alpha-1}-s^{-\alpha-1}|\,ds,
    \end{align*}
    and
    \[
        \int_0^t s^{\lambda - \alpha - 1} \, ds \leq t^{\lambda - \alpha} \leq
        \begin{cases}
           h^{\lambda - \alpha} \leq T h^{\lambda - \alpha-1} & h \geq t,\\
              T^\lambda h^{-\alpha} \leq T^\lambda h^{\lambda - \alpha-1} & h < t.
        \end{cases}
    \]

\subsubsection*{Proof of the continuity of $\mathcal{I}_2$ from $C^\lambda_0$ to $C^{\lambda-\alpha}_0$}
    Let $h\in(0,1)$ and $t\in(0,T)$ such that $t+h\leq T$. 
    We have
    \[
    |\mathcal{I}_2g(t)|\leq \|g\|_{C_{0}^{\lambda}}\sup_{s\in(0,T]}t^{\alpha}|d_2(s)|\,\int_0^{t}s^{\lambda-\alpha}\,ds,
    \]
    so $\mathcal{I}_2f$ can be defined on $[0,T]$ and $\mathcal{I}_2 g(0)=0$.
    Next we evaluate the difference. We have
        \begin{align*}
            |&\mathcal{I}_2g(t+h)-\mathcal{I}_2g(t)|\\
            &=|(g(t+h)-g(t))\int_{0}^{t+h}d_2(s)\,ds+g(t)\int_{t}^{t+h}d_2(s)\,ds\\
            &\quad-\int_{0}^{t}(g(t+h-s)-g(t-s))d_2(s)\,ds-\int_{t}^{t+h}g(t+h-s)d_2(s)\,ds|\\
            &\leq \|g\|_{C_{0}^{\lambda}}(
            h^{\lambda}\int_{0}^{t+h}|d_2(s)|\,ds+\int_{t}^{t+h}|d_2(s)|\,ds\\
            &\quad\quad\quad\quad +h^{\lambda}\int_{0}^{t+h}|d_2(s)|\,ds+\int_{t}^{t+h}|d_2(s)|\,ds),
        \end{align*}
    which is bounded by a constant multiple of $\|g\|_{C_0^{\lambda}}h^{\lambda-\alpha}$ from the following inequalities:
    \[
        \int_{0}^{t+h}|d_2(s)|\,ds\leq C\int_{0}^{t+h} s^{-\alpha}\,ds\leq C(t+h)^{1-\alpha}\leq C T h^{-\alpha}
    \]
    \[
        \int_{t}^{t+h}|d_2(s)|\,ds
        \leq C\int_{t}^{t+h}s^{-\alpha}\,ds
        \leq C((t+h)^{1-\alpha}-t^{1-\alpha})\,ds
        \leq C h^{1-\alpha}.
    \]
    Here we have used that $(\cdot)^{1-\alpha}$ is H\"older continuous.
\hfill{$\square$}

\end{document}